\documentclass{amsart}
\newcommand*{\mailto}[1]{\href{mailto:#1}{\nolinkurl{#1}}}
\usepackage{hyperref}
\newcommand{\arxiv}[1]{\href{http://arxiv.org/abs/#1}{arXiv:#1}}

\newtheorem{theorem}{Theorem}[section]

\newtheorem{lemma}[theorem]{Lemma}
\newtheorem{corollary}[theorem]{Corollary}
\newtheorem{remark}[theorem]{Remark}
\newtheorem{hypothesis}[theorem]{Hypothesis}

\newcommand{\R}{{\mathbb R}}
\newcommand{\N}{{\mathbb N}}
\newcommand{\Z}{{\mathbb Z}}
\newcommand{\C}{{\mathbb C}}

\newcommand{\nn}{\nonumber}
\newcommand{\be}{\begin{equation}}
\newcommand{\ee}{\end{equation}}

\newcommand{\ti}{\tilde}

\newcommand{\E}{\mathrm{e}}
\newcommand{\I}{\mathrm{i}}

\newcommand{\im}{\mathrm{Im}}

\newcommand{\eps}{\varepsilon}

\newcommand{\lam}{\lambda}
\newcommand{\gam}{\gamma}

\newcommand{\ta}{\theta}

\numberwithin{equation}{section}

\begin{document}

\title[Inverse Eigenvalue Problems for Spherical Schr\"odinger Operators]{Inverse Eigenvalue Problems for Perturbed Spherical Schr\"odinger Operators}

\author[A. Kostenko]{Aleksey Kostenko}
\address{Institute of Applied Mathematics and Mechanics\\
NAS of Ukraine\\ R. Luxemburg str. 74\\
Donetsk 83114\\ Ukraine\\ and School of Mathematical Sciences\\
Dublin Institute of Technology\\
Kevin Street\\ Dublin 8\\ Ireland}
\email{\mailto{duzer80@gmail.com}}

\author[A. Sakhnovich]{Alexander Sakhnovich}
\address{Faculty of Mathematics\\
Nordbergstrasse 15\\ 1090 Wien\\ Austria}
\email{\mailto{Oleksandr.Sakhnovych@univie.ac.at}}
\urladdr{\url{http://www.mat.univie.ac.at/~sakhnov/}}

\author[G.\ Teschl]{Gerald Teschl}
\address{Faculty of Mathematics\\
Nordbergstrasse 15\\ 1090 Wien\\ Austria\\ and International
Erwin Schr\"odinger
Institute for Mathematical Physics\\ Boltzmanngasse 9\\ 1090 Wien\\ Austria}
\email{\mailto{Gerald.Teschl@univie.ac.at}}
\urladdr{\url{http://www.mat.univie.ac.at/~gerald/}}

\thanks{Inverse Problems {\bf 26}, 105013, 14pp (2010)}
\thanks{{\it Research supported by the Austrian Science Fund (FWF) under Grant No.\ Y330}}

\keywords{Schr\"odinger operators, spectral theory, strongly singular potentials}
\subjclass[2000]{Primary 34B20, 34L15; Secondary 81V45 , 47A10}

\begin{abstract}
We investigate the eigenvalues of perturbed spherical Schr\"odinger operators under the
assumption that the perturbation $q(x)$ satisfies $x q(x) \in L^1(0,1)$. We show that the
square roots of eigenvalues are given by the square roots of the unperturbed eigenvalues
up to an decaying error depending on the behavior of $q(x)$ near $x=0$. Furthermore,
we provide sets of spectral data which uniquely determine $q(x)$.
\end{abstract}

\maketitle

\section{Introduction}
\label{sec:int}

Given a Schr\"odinger operator with a rotationally symmetric potential, separation of variables leads
to the spherical Schr\"odinger operator (e.g., \cite{tschroe}, \cite{wdln})
\be
H = - \frac{d^2}{dx^2} + \frac{l(l+1)}{x^2} +q(x), \qquad l=0,1,2,\dots.
\ee
In this note we are interested in the case where the particle is confined to a finite ball (of radius $1$
for notational simplicity --- which can always be achieved by scaling $x$). This problem has attracted much
interest in the past and several results concerning the eigenvalues of these problems have been derived.
The first results are based on the seminal work by Guillot and Ralston \cite{gr} which deals with the case $l=1$
and $q \in L^2(0,1)$. It was later extended by Carlson \cite{car}, \cite{car2} and recently by Serier \cite{se} who
extended their results to arbitrary $l\in\N_0$. However, the assumption $q \in L^2(0,1)$ clearly excludes the physically
interesting case of a Coulomb type singularity $q(x) = \frac{\gam}{x} + \dots$. This case was included for $l=0$
in the work of Savchuk and Shkalikov \cite{ss}, \cite{ss2} who considered $q\in W^{-1,2}(0,1)$. Their work
was later extended by Albeverio, Hryniv, and Mykytyuk who first covered the case $q\in W^{-1,p}(0,1)$, $p\in[1,\infty)$, for $l=0$
in \cite{ahm} and later on extended this to all $l\in\N_0$ in \cite{ahm2} using the double commutation method. Their
condition includes in particular the case $x q(x) \in L^p(0,1)$, $p\in[1,\infty)$, which will be the condition we are interested in here.

Our main motivation for this paper is the paper by Zhornitskaya and Serov \cite{zs} who treat the general case $l \ge -\frac{1}{2}$
under the assumption $q \in L^1(0,1)$. We want to extend their results in
several ways: First of all we will replace the condition $q(x) \in L^1(0,1)$ by $x q(x) \in L^1(0,1)$.
Moreover, they show in \cite{zs} that the Dirichlet eigenvalues satisfy
\be
\mu_n = \left(j_{l+\frac{1}{2},n} + \eps_n \right)^2
\ee
where $j_{l+\frac{1}{2},n} =  \pi (n +\frac{l}{2})+ O(n^{-1})$ are the zeros of the Bessel function $J_{l+1/2}(z)$ and the error
satisfies $|\eps_n| \leq \frac{\pi}{4}$ (this is claimed for all $n$ but only proven for large $n$). We will show that
the error satisfies
\be
\eps_n = O\left(\int_0^1 \frac{y |q(y)|}{1 + n y} dy\right).
\ee
(For $l=-\frac{1}{2}$ one has to replace $q(y)$ by $(1-\log(y))q(y)$.)
In particular, for $q\in L^1(0,1)$ we get $\eps_n = O(n^{-1})$ and for a Coulomb type singularity $q(x) = \frac{\gam}{x} + L^1(0,1)$
we get $\eps_n = O(n^{-1} \log(n))$.

Based on this information we will give some sets of spectral data which, in addition to the Dirichlet spectrum, uniquely determine
$q$ and $l$ again generalizing the corresponding result from \cite{zs} to the case $x q(x) \in L^1(0,1)$. Moreover, in their construction
they use the fact that a certain Wronskian does not vanish. Unfortunately this Wronskian can indeed vanish (we will give a simple
counter example due to Ralston in Remark~\ref{rem:wro}) but we will show that the use of this {\em fact} can be avoided.

\section{The spherical Schr\"odinger operator}

Our prototypical example will be the spherical Schr\"odinger operator given by
\be
H_l = - \frac{d^2}{dx^2} + \frac{l(l+1)}{x^2}, \qquad x\in (0,1), \quad l \ge -\frac{1}{2}.
\ee
Note that we explicitly allow non-integer values of $l$ such that we also cover the case of
arbitrary space dimension $n\ge 2$, where $l(l+1)$ has to be replaced by $l(l+n-2) + (n-1)(n-3)/4$ \cite[Sec.~17.F]{wdln}.

With the usual boundary conditions at $x=0$ (for $l\in[-\frac{1}{2},\frac{1}{2})$) and $x=1$
\be
\lim_{x\to0} x^l ( (l+1)f(x) - x f'(x))=0, \qquad f(1) =0 \text{ or } f'(1) + \beta f(1) =0,
\ee
it gives rise to a self-adjoint operator in the Hilbert space $L^2(0,1)$.
Two linearly independent solutions of
\be
- f''(x) + \frac{l(l+1)}{x^2} f(x) = z f(x)
\ee
are given by (see \cite[(9.1.49)]{as})
\be\label{defphil}
\phi_l(z,x) = z^{-\frac{2l+1}{4}} \sqrt{\frac{\pi x}{2}} J_{l+\frac{1}{2}}(\sqrt{z} x),
\ee
\be\label{defthetal}
\theta_l(z,x) = -z^{\frac{2l+1}{4}} \sqrt{\frac{\pi x}{2}} \begin{cases}
\frac{-1}{\sin((l+\frac{1}{2})\pi)} J_{-l-\frac{1}{2}}(\sqrt{z} x), & {l+\frac{1}{2}} \in \R_+\setminus \N_0,\\
Y_{l+\frac{1}{2}}(\sqrt{z} x) -\frac{1}{\pi}\log(z) J_{l+\frac{1}{2}}(\sqrt{z} x), & {l+\frac{1}{2}} \in\N_0,\end{cases}
\ee
where $J_{l+1/2}$ and $Y_{l+1/2}$ are the usual Bessel and Neumann functions.
All branch cuts are chosen along the negative real axis unless explicitly stated otherwise.
If $l$ is an integer they of course reduce to spherical Bessel and Neumann functions
and can be expressed in terms of trigonometric functions (cf.\ e.g.\ \cite{as, wa} and also \cite[Sect.~10.4]{tschroe}).

Using the power series for the Bessel and Neumann functions one verifies that they have the form
\be
\phi_l(z,x) =  x^{l+1} \frac{\sqrt{\pi}}{\Gamma(l+\frac{3}{2}) 2^{l+1}} f_l(z x^2),
\ee
\be
\theta_l(z,x) =  \frac{\Gamma(l+\frac{3}{2}) 2^{l+1}}{x^{l}\sqrt{\pi}} \begin{cases}
\frac{1}{2l+1} g_l(z x^2), & l+\frac{1}{2} \in \R_+\setminus \N_0,\\
\frac{1}{2l+1} g_l(z x^2) -\frac{ (z x^2)^{l+\frac{1}{2}}\log(x)}{\Gamma(l+\frac{3}{2})^2 2^{2l+1}}  f_l(z x^2), & l+\frac{1}{2} \in\N,\\
(\log(2)-\gam) g_l(z x^2) - \log(x) f_l(z x^2), & l =-\frac{1}{2},\end{cases}
\ee
where $f_l(z)$, $g_l(z)$ are entire functions with $f_l(0)=g_l(0)=1$ and $\gam$ is the Euler--Mascheroni constant.

In particular,  both functions are entire and according to \cite[(9.1.16)]{as} their Wronskian is given by
\be \label{a0}
W(\theta_l(z),\phi_l(z))=1.
\ee
The eigenvalues of $H_l$ with a Dirichlet boundary condition at $x=1$ are given by the
zeros of the entire function $\phi_l(z,1)$ which are the squares of the positive zeros of the Bessel function
$J_{l+1/2}$ of order $l+1/2$:
\be
\mu_{l,n} = (j_{l+1/2,n})^2.
\ee
Similarly, the eigenvalues of $H_l$ with the boundary condition $f'(1)+\beta f(1)=0$
are given by the zeros of the entire function $\phi_l'(z,1)+ \beta \phi_l(z,1)$ and the positive eigenvalues
are the squares of the positive zeros of $zJ_{l+3/2}(z) - (\beta+l+1) J_{l+1/2}(z)$ (cf. \cite[\S 3.2]{wa}):
\be
\lambda^\beta_{l,n} = (j^\beta_{l+1/2,n})^2.
\ee
Here the eigenvalues are counted according to $\lambda^\beta_{l,0} < \mu_{l,1} < \lambda^\beta_{l,1} < \dots$.
In particular, observe that the first eigenvalue $\lambda^\beta_{l,0}$ will be zero for $\beta =-(l+1)$ and negative for
$\beta <-(l+1)$.

Note that asymptotically (\cite[(9.5.12)]{as})
\begin{align}\label{a6}
\sqrt{\mu_{l,n}}&=j_{l+\frac{1}{2},n} = \left(n +\frac{l}{2}\right) \pi + O(n^{-1}), \qquad\\\label{a7}
\sqrt{\lam^\beta_{l,n}}&=j^\beta_{l+\frac{1}{2},n} =j_{l+\frac{3}{2},n}+O(n^{-1})= \left(n +\frac{l+1}{2} \right) \pi + O(n^{-1})
\end{align}
for fixed $l$.

Now let us look at perturbations
\be
H = H_l + q(x)
\ee
assuming that the potential $q$ satisfies the following conditions:

\begin{hypothesis}\label{hyp_q}
Let $l\in [-\frac{1}{2},\infty)$. Set
\be \label{f1}
\ti{q}(x) = \begin{cases}
|q(x)|, & l > -\frac{1}{2},\\
(1-\log(x)) |q(x)|, & l = -\frac{1}{2},
\end{cases}
\ee
and suppose $q$ is real-valued such that
\be\label{a0!}
\begin{array}{c}
x \ti{q}(x) \in L^1(0,1).
\end{array}
\ee
\end{hypothesis}

\begin{lemma}\label{lemphi}
Assume Hypothesis \ref{hyp_q}.
Then there is a solution $\phi(z,x)$ of $H f = z f$ which is entire with respect to $z$ and satisfies the integral
equation
\be \label{a1}
\phi(z,x) = \phi_l(z,x) + \int_0^x G_l(z,x,y) q(y) \phi(z,y) dy,
\ee
where
\be \label{a2}
G_l(z,x,y) = \phi_l(z,x) \theta_l(z,y) - \phi_l(z,y) \theta_l(z,x)
\ee
is the Green function of the initial value problem. Moreover, this solution satisfies the estimate
\be\label{estphi}
| \phi(z,x) - \phi_l(z,x)| \leq C \left(\frac{x}{1+ |z|^{1/2} x}\right)^{l+1} \E^{|\im(z^{1/2})| x} \int_0^x \frac{y \ti{q}(y)}{1 +|z|^{1/2} y} dy.
\ee
The derivative is given by
\be  \label{a3}
\phi'(z,x) = \phi_l'(z,x) + \int_0^x \frac{\partial}{\partial x}G_l(z,x,y) q(y) \phi(z,y) dy
\ee
and satisfies the estimate
\be\label{estphi'}
| \phi'(z,x) - \phi_l'(z,x)| \leq C  \left(\frac{x}{1+ |z|^{1/2} x}\right)^l \E^{|\im(z^{1/2})| x} \int_0^x \frac{y \ti{q}(y)}{1 +|z|^{1/2} y} dy.
\ee
\end{lemma}

\begin{proof}
In a way similar to \cite{gr} (see also \cite{se})
this can be shown by iteration using Lemmas \ref{LaA1} and \ref{LaA2}.
Namely, it is easy to show that $\phi$ given by
\begin{align}&\label{iter1}
 \phi=\sum_{n=0}^{\infty}\phi_{l,n}, \quad \phi_{l,0}:=\phi_l, \\
&\label{iter2}
\phi_{l,n+1}(z,x):=\int_0^x G_l(z,x,y) q(y) \phi_{l,n}(z,y) dy, \quad n \in \N_0,
\end{align}
satisfies \eqref{a1}. The inequalities
\begin{align}&\label{iter3}
|\phi_{l,n}(z,x)|\leq \frac{C^{n+1}}{n!}\left(\frac{x}{1+ |z|^{1/2} x}\right)^{l+1} \E^{|\im(z^{1/2})| x}
\left(\int_0^x \frac{y \ti{q}(y)}{1 +|z|^{1/2} y} dy\right)^n,
\end{align}
which are necessary to prove the convergence in \eqref{iter1}, follow by induction from  \eqref{a0!} and Lemma \ref{LaA1}.
By \eqref{iter1} and \eqref{iter3}, we get inequality \eqref{estphi}, where $\ti{q}$ is given in \eqref{f1}.
Hear $C$ is to be understood as a generic constant whose value is different
in \eqref{estphi} and in Lemma \ref{LaA1}. The properties of $\phi^{\prime}$ follow in a quite
similar way from \eqref{iter1}--\eqref{iter3} and Lemma \ref{LaA2}.
Finally, using \eqref{a0}, \eqref{a1}, \eqref{a2}, and \eqref{a3} one can see that $H\phi=z\phi$.
\end{proof}

Using the well-known asymptotic formulas for Bessel function \cite[(9.2.1), (9.2.11)]{as} we obtain
\begin{align}\label{asymphibes}
\phi(z,x) &= z^{-\frac{l+1}{2}}\left(\sin\bigl(\sqrt{z} x- \frac{l \pi}{2}\bigr) + O\bigl(|z|^{-1/2}\E^{x |\im(\sqrt{z})|}\bigr)\right),\\\label{asymphibes2}
\phi'(z,x) &= z^{-\frac{l}{2}}\left(\cos\bigl(\sqrt{z} x- \frac{l \pi}{2}\bigr) + O\bigl(|z|^{-1/2}\E^{x |\im(\sqrt{z})|}\bigr)\right),
\end{align}
as $z\to\infty$. Next, note that
\be\label{limphi}
\lim_{x\to 0}  x^{-l-1} \phi(z,x) = \lim_{x\to 0}  x^{-l-1} \phi_l(z,x) = \frac{\sqrt{\pi}}{\Gamma(l+\frac{3}{2}) 2^{l+1}}.
\ee
Moreover, using  d'Alembert's formula (cf.\ \cite[Sect.~XI.6]{har}) a second linearly independent solution,
satisfying $W(\theta(z),\phi(z))=1$, is given by
\be
\theta(z,x) = - \phi(z,x) \int_x^c \frac{dy}{\phi(z,y)^2}.
\ee
where $c=c(z)$ has to be chosen such that $\phi(z,x)$ does not vanish in $(0,c)$. In particular it is
straightforward to show

\begin{corollary}
Assume Hypothesis \ref{hyp_q}.
The differential equation $H f = z f$ has two linearly independent solutions,
satisfying $W(\theta(z),\phi(z))=1$, of the form
\be
\phi(z,x) = x^{l+1} \ti{\phi}(z,x), \qquad \theta(z,x) = \begin{cases}
\frac{x^{-l}}{2l+1} \ti{\theta}(z,x), & l > -\frac{1}{2},\\
-x^{1/2} \log(x) \ti{\theta}(z,x), & l = -\frac{1}{2},
\end{cases}
\ee
where $\ti{\phi}(z,x), \ti{\theta}(z,x) \in C(\C,[0,1])$ are jointly continuous and $\ti{\phi}(z,0) =\ti{\theta}(z,0)^{-1} \ne 0$.
\end{corollary}

Unfortunately, since $c(z)\to 0$ as $z\to\infty$, this simple approach to get a second solution
looses control over $\theta(z,x)$ as a function of $z$. In particular, it is not clear that this
second solution can be chosen to be entire as a function of $z$. We will construct a second solution
with better control with respect to $z$ in Lemma~\ref{lempsi} below.

In any case, the behavior of the solutions of $H f = z f$ near $x=0$ implies

\begin{theorem}
Assume Hypothesis \ref{hyp_q}. The differential equation $H = H_l +q$ is limit circle at $x=0$ if $l\in[-\frac{1}{2},\frac{1}{2})$ and
limit point at $x=0$ for $l\geq \frac{1}{2}$.  In particular, $H$ associated with the boundary conditions at $x=0$
(for $l\in[-\frac{1}{2},\frac{1}{2})$) and $x=1$
\be
\lim_{x\to0} x^l ( (l+1)f(x) - x f'(x))=0, \qquad f(1) =0 \text{ or } f'(1) + \beta f(1) =0
\ee
is self-adjoint. Moreover, the spectrum of $H$ is purely discrete and bounded from below.
\end{theorem}

\begin{proof}
Since $\phi(z,x)$ and $\theta(z,x)$ are both in $L^2(0,1)$ if and only if $l\in[-\frac{1}{2},\frac{1}{2})$ we deduce that $H$ is limit circle at $x$
if and only if $l\in[-\frac{1}{2},\frac{1}{2})$. Moreover, in this case we can choose the boundary condition (cf.\ \cite[Sect~9.2]{tschroe})
\[
\lim_{x\to0} W(\phi(0,x),f(x)) = \frac{\sqrt{\pi}}{\Gamma(l+\frac{3}{2}) 2^{l+1}} \lim_{x\to0} x^l ( (l+1)f(x) - x f'(x))=0,
\]
where we have used \eqref{limphi} and the fact that any solution in the maximal domain of the differential
expression satisfies $\lim_{x\to 0} x^{l+1} f(x) =0$.

Furthermore, since $\phi(z,x)$ has only a finite number of zeros inside $(0,1)$,
the differential expression is nonoscillatory for every $z\in\R$, we conclude that the
spectrum is purely discrete and bounded from below (cf.\ \cite[Thm.~14.9]{wdln}).
\end{proof}

Now we are able to prove our eigenvalue asymptotics using a refined version of the approach by
P\"oschel and Trubowitz \cite{pt}.

\begin{theorem} \label{thmev}
Assume Hypothesis \ref{hyp_q}.
Then the Dirichlet eigenvalues corresponding to the boundary condition $f(1)=0$ satisfy
\be
\mu_n = \left(j_{l+\frac{1}{2},n} + \eps_n \right)^2
\ee
where
\[
\eps_n = O\left(\int_0^1 \frac{y \ti{q}(y)}{1 + n y} dy\right).
\]
Similarly the eigenvalues corresponding to the boundary condition $f'(1)+\beta f(1)=0$ satisfy
\be
\lam_n^\beta = \left(j^\beta_{l+\frac{1}{2},n}  + \tilde\eps_n \right)^2
\ee
where $\tilde\eps_n$ is of the same order as $\eps_n$.
\end{theorem}

\begin{proof}
We set $\phi_l(z):= \phi_l(z,1) = \sqrt{\pi /2}\ z^{-\frac{2l+1}{4}} J_{l+1/2}(\sqrt{z})$ and $\phi(z)= \phi(z,1)$.
Then our estimate \eqref{estphi} reads
\[
|\phi(z) - \phi_l(z)| \leq C \frac{\E^{|\im \sqrt{z}|}}{(1+\sqrt{|z|})^{l+1}} \eps(z),\quad \eps(z):=\int_0^1 \frac{y \ti{q}(y)}{1 +|z|^{1/2} y} dy.
\]
 Next, using \cite[(9.2.1), (9.2.11)]{as}, we have
\begin{align}
\phi_l(z) &= z^{-\frac{l+1}{2}}\left(\sin\bigl(\sqrt{z} - \frac{l \pi}{2}\bigr) + O\bigl(|z|^{-1/2}\E^{|\im(z^{1/2})|}\bigr)\right), \label{estphi_l} \\\label{estphi'_l}
\dot{\phi}_l(z) &= \frac{1}{2}z^{-\frac{l+2}{2}}\left(\cos\bigl(\sqrt{z} - \frac{l \pi}{2}\bigr) + O\bigl(|z|^{-1/2}\E^{|\im(z^{1/2})|}\bigr)\right),
\end{align}
where the dot denotes a derivative with respect to $z$.
Hence, taking into account \eqref{a6}, \eqref{estphi_l}, and \eqref{estphi'_l}, we get
\be\label{estimates}
\phi_l\bigl((j_{l+\frac{1}{2},n})^2\bigr) = 0, \quad \dot{\phi}_l\bigl((j_{l+\frac{1}{2},n})^2\bigr) = \frac{1}{2}(j_{l+\frac{1}{2},n})^{-(l+2)}\left((-1)^n + O(n^{-1})\right).
\ee
Furthermore, \eqref{estphi_l}--\eqref{estimates} together with the mean value theorem yield (for sufficiently large $n$)
\begin{align*}
\big|\phi_l\bigl((j_{l+\frac{1}{2},n}\pm \eps_n)^2\bigr)\big| \geq \frac{1}{4} (j_{l+\frac{1}{2},n}+\eps_n)^{-(l+2)}\big|(j_{l+\frac{1}{2},n}\pm \eps_n)^2-j_{l+\frac{1}{2},n}^2\big|\\
\ge \frac{1}{4}\frac{\eps_n|2j_{l+\frac{1}{2},n}-\eps_n|}{(j_{l+\frac{1}{2},n}+\eps_n)^{l+2}} > \frac{1}{4}\frac{\eps_n}{(j_{l+\frac{1}{2},n}+1)^{l+1}}
=2C\frac{\eps(j_{l+\frac{1}{2},n})}{(j_{l+\frac{1}{2},n}+1)^{l+1}},
\end{align*}
where $\eps_n= 8 C \eps(j_{l+\frac{1}{2},n})$. Thus
\begin{align}
\big|\phi\bigl((j_{l+\frac{1}{2},n}\pm \eps_n)^2\bigr) - \phi_l\bigl((j_{l+\frac{1}{2},n}\pm \eps_n)^2\bigr)\big| \leq&
C \frac{\eps(j_{l+\frac{1}{2},n}-\eps_n)}{(1+j_{l+\frac{1}{2},n}-\eps_n)^{l+1}} \nn\\ &< \big|\phi_l\bigl((j_{l+\frac{1}{2},n}\pm \eps_n)^2\bigr)\big|.\nn
\end{align}
This shows that $\phi(z)$ has different signs at $(j_{l+\frac{1}{2},n} - \eps_n)^2$ and $(j_{l+\frac{1}{2},n} + \eps_n)^2$, and thus
there is at least one zero in between.

Let us show that $\phi$ has no other zeros. Since $|\sin z|>\frac{1}{4}\E^{|\im z|}$ if $|z-\pi n|\geq \frac{\pi}{4}$ for all $n\in\Z$
(see for instance \cite[Lemma 2.1]{pt}), there exists a $N\in \N$ such that
\[
|\phi_l(z)|\geq C_l\frac{\E^{|\im \sqrt{z}|}}{|z|^{\frac{l+1}{2}}},\qquad |z|\geq N,\quad \bigl|\sqrt{z}-\pi n -\frac{l\pi}{2}\bigr|\geq \frac{\pi}{4},\quad n\ge N,
\]
with some positive constant $C_l>0$ independent of $z$ (this estimate can also be deduced from \cite[Lemma 22.1]{lev}).
Since $\eps(z)\to0$ as $|z|\to\infty$, there exists $K>0$ such that $\eps(z)<(C_l C)^{-1}$ if $|z|>K$.
Thus, on contours $|z|=(K+N+\frac{l+1}{2})\pi$ and $\bigl|\sqrt{z}-\pi n -\frac{l\pi}{2}\bigr|= \frac{\pi}{4}$, $n\geq K+N+1$, we obtain
\[
| \phi(z) - \phi_l(z)|<|\phi_l(z)|.
\]
By Rouch\'e's theorem, $\phi(z)$ has as many roots as $\phi_l(z)$ in each of the bounded regions and the remaining unbounded region.
Since all roots are simple, we are done.

An analogous argument can be given for the general eigenvalues based on zeros of
\[
\phi_l^\beta(z) = \phi_l'(z,1) + \beta \phi_l(z,1) =
\sqrt{\frac{\pi}{2}} z^{-\frac{2l+1}{4}} \left((\beta+l+1) J_{l+\frac{1}{2}}(\sqrt{z}) - \sqrt{z}J_{l+\frac{3}{2}}(\sqrt{z}) \right).
\]
\end{proof}

Similar results for $q\in L^2$ were given by Guillot and Ralston \cite{gr}, Carlson \cite{car}, Serier \cite{se}.
Note that the case $x q(x) \in L^p$ for $l\in\N$ is covered in \cite{ahm, ss, ss2} for $l=0$ and in
\cite{ahm2} for $l\in\N_0$, where the direct and inverse spectral problems have been effectively studied for
operators with distributional potentials $q\in W^{-1,p}(0,1)$.

Finally we come to uniqueness results for the inverse problem. Considering the solution
\be\label{a4i}
\psi_l^\infty(z,x) = \theta_l(z,1) \phi_l(z,x) - \phi_l(z,1) \theta_l(z,x) = G_l(z,x,1)
\ee
satisfying the initial conditions $(\psi_l^\infty(z,1),(\psi_l^\infty(z,1))')=(0,1)$ and
\begin{align}\nn
\psi_l^\beta(z,x) &= - \big(\beta \theta_l(z,1)+ \theta_l'(z,1)\big) \phi_l(z,x) + \big(\beta \phi_l(z,1)+ \phi_l'(z,1)\big) \theta_l(z,x)\\ \label{a4}
&= - \beta G_l(z,x,1) + \frac{\partial G_l}{\partial y}(z,x,1)
\end{align}
satisfying the initial conditions $(\psi_l^\beta(z,1),(\psi_l^\beta(z,1))')=(1,-\beta)$, we obtain the analog of Lemma~\ref{lemphi}.

\begin{lemma}\label{lempsi}
Assume Hypothesis \ref{hyp_q}.
Then there is a solution $\psi^\beta(z,x)$ of $H u = z u$ which is entire with respect to $z$ and satisfies the integral
equation
\be \label{a5}
\psi^\beta(z,x) = \psi_l^\beta(z,x) - \int_x^1 G_l(z,x,y) q(y) \psi^\beta(z,y) dy.
\ee
Moreover, for $l>-1/2$ this solution satisfies the estimate
\begin{align}\label{estpsi_B}
| \psi^\infty(z,x) - \psi_l^\infty(z,x)| \leq C \left(\frac{1+ |z|^{1/2} x}{x+ |z|^{1/2}x}\right)^{l}
 \frac{\E^{|\im(z^{1/2})| (1-x)}}{1+|z|^{1/2}} \int_x^1 \frac{y \ti{q}(y)}{1 +|z|^{1/2} y} dy,\\ \label{estpsi}
| \psi^\beta(z,x) - \psi_l^\beta(z,x)| \leq C \left(\frac{1+ |z|^{1/2} x}{x+ |z|^{1/2}x}\right)^{l}
 \E^{|\im(z^{1/2})| (1-x)} \int_x^1 \frac{y \ti{q}(y)}{1 +|z|^{1/2} y} dy.
\end{align}
In the case $l=-\frac{1}{2}$ an additional factor $ (1-\log(x))$ appears in the right-hand sides of
\eqref{estpsi_B} and  \eqref{estpsi}.
\end{lemma}

\begin{proof}
Suppose $\beta\neq \infty$ first.
In a way quite similar to the proof of Lemma \ref{lemphi} one can show that
the solution $\psi^{\beta}$ admits representation
\begin{align}&\label{iter1p}
\psi^{\beta}=\sum_{n=0}^{\infty}\psi^{\beta}_{l,n}, \quad \psi^{\beta}_{l,0}:=\psi^{\beta}_l, \\
&\label{iter2p}
\psi^{\beta}_{l,n+1}(z,x):=-\int_x^1 G_l(z,x,y) q(y) \psi^{\beta}_{l,n}(z,y) dy \quad (n \in \N_0),
\end{align}
where the functions $\psi^{\beta}_{l,n}$ satisfy inequalities
\begin{align}&\label{iter3p}
|\psi^{\beta}_{l,n}(z,x)|\leq \frac{C^{n+1}}{n!} \left(\frac{1+ |z|^{1/2} x}{x+ |z|^{1/2}x}\right)^{l} \E^{|\im(z^{1/2})| (1-x)}
\left(\int_x^1 \frac{y \ti{q}(y)}{1 +|z|^{1/2} y} dy\right)^n
\end{align}
for $l>-\frac{1}{2}$ and inequalities
\begin{align}\label{iter4p}
|\psi^{\beta}_{l,n}(z,x)|\leq \frac{C^{n+1}}{n!}  (1-\log(x)) \left(\frac{x+ |z|^{1/2} x}{1+ |z|^{1/2}x}\right)^{1/2}\\
\times\
\E^{|\im(z^{1/2})| (1-x)}
\left(\int_x^1 \frac{y \ti{q}(y)}{1 +|z|^{1/2} y} dy\right)^n\nn
\end{align}
for $l=-\frac{1}{2}$.
Indeed, the inequality \eqref{iter3p} for $n=0$ and for some $C>0$ follows from \eqref{a4}, \eqref{estGlp},
and the fact that according to \eqref{a2} we have $G_l(z,x,y)= - G_l(z,y,x)$.
Next, we use \eqref{iter2p} and \eqref{estGl} to prove \eqref{iter3p} for all $n$ by induction.
By  \eqref{iter1p}--\eqref{iter3p}
we get  \eqref{a5} and  \eqref{estpsi}.
Finally, according to \eqref{a5}  we have $H\psi^{\beta}=z\psi^{\beta}$ and
$(\psi^\beta(z,1),(\psi^\beta(z,1))')=(\psi_l^\beta(z,1),(\psi_l^\beta(z,1))')=(1,-\beta)$.

The inequality \eqref{iter4p} for $n=0$ follows from \eqref{a4} and  \eqref{a21},
and the remaining part  of  the proof for $l=-1/2$ is analogous to the case $l>-1/2$.

The case $\beta=\infty$ can be treated in a similar way, it suffices to note that
\be\label{iter3pi}
|\psi^{\infty}_{l,0}(z,x)|\leq \frac{C}{1+ |z|^{1/2}} \left(\frac{1+ |z|^{1/2} x}{x+ |z|^{1/2}x}\right)^{l} \E^{|\im(z^{1/2})| (1-x)}
\ee
for $l>-\frac{1}{2}$ and
\be\label{iter4pi}
|\psi^{\infty}_{-1/2,0}(z,x)|\leq  (1-\log(x))\frac{C}{1+ |z|^{1/2}} \left(\frac{x+ |z|^{1/2}x}{1+ |z|^{1/2}x}\right)^{1/2}
\E^{|\im(z^{1/2})| (1-x)}
\ee
for $l=\frac{1}{2}$.
\end{proof}

\begin{remark}\label{rem:wro}
In \cite{zs} a second solution $\tilde{\psi}(z,x)$ is constructed by considering
\be
\tilde{\psi}(z,x) = \theta_l(z,x) - \int_x^1 G_l(z,x,y) q(y) \tilde{\psi}(z,y) dy.
\ee
It is claimed to be linearly independent for all $z$ and a reference is made to \cite{gr}, where
the corresponding claim for $l=1$ was made. However, this is wrong as the following
counter example shows (a similar counter example was communicated to us by Ralston; see also \cite{gr2}):

The function $\tilde{\psi}$ satisfies $(\tilde{\psi}(z,1),\tilde{\psi}'(z,1))= (\theta_l(z,1),\theta_l'(z,1))$.
Let $y_{l+1/2,n}$ be a zero of $\theta_l(1,z)$ and chose $q(x)$ such that $y_{l+1/2,n}$ is an Dirichlet eigenvalue
of $H$. Then $W(\tilde{\psi}(y_{l+1/2,n}), \phi(y_{l+1/2,n}))=0$.
\end{remark}

Now we come to our uniqueness results:
\begin{theorem}
Assume Hypothesis \ref{hyp_q}.
The following set of spectral data determine $q$ and $l$ uniquely:
\begin{enumerate}
\item
Two sets of eigenvalues $\lam^\alpha_n$ and $\lam^\beta_n$ for $\alpha\neq \beta$ (including the case $\lam^\infty_n = \mu_n$).
\item
The Dirichlet eigenvalues $\{ \mu_n \}_{n=1}^\infty$ together with the norming constants
\[
\gam_n^{-1} = \frac{1}{ \phi'(\mu_n,1)^2} \int_0^1 \phi(\mu_n,x)^2 dx = -\frac{\dot{\phi}(\mu_n,1)}{\phi'(\mu_n,1)},
\]
or the set of eigenvalues $\{\lam_n^\beta\}_{n=0}^\infty$ together with the norming constants
\[
(\gam_n^\beta)^{-1} = \frac{1}{ \phi^\beta(\lam_n^\beta)^2} \int_0^1 \phi(\lam_n^\beta,x)^2 dx = \frac{\dot{\phi^\beta}(\lam_n^\beta)}{(1+\beta^2)\phi(\lam_n^\beta,1)},
\]
where $\phi^\beta(z):=\phi'(z,1)+\beta\phi(z,1)$.
\item
The eigenvalues $\lam^\beta_n$ together with $\phi(\lam_n^\beta,1)$ if $\beta\ne\infty$ or ${\phi}'(\lam_n^\beta,1)$ if $\beta\ne 0$.
\end{enumerate}
Here the dot and prime denote derivatives with respect to $z$ and $x$, respectively.
\end{theorem}

\begin{proof}
(i) To see the first claim recall the Weyl $m$-function of $H$ at the regular endpoint $x=1$, which is given by
\[
m_\beta(z) = \frac{\phi(z,1)- \beta \phi'(z,1)}{\phi'(z,1) + \beta \phi(z,1)}
\]
and consider
\[
\ti{m}(z) = m_\beta(z) - \frac{\alpha\beta +1}{\beta - \alpha} = \frac{1+\beta^2}{\alpha-\beta}\cdot\frac{\phi(z,1) +\alpha \phi'(z,1)}{\phi'(z,1) + \beta \phi(z,1)}.
\]
Then the well-known asymptotics \cite{tschroe}
\[
m_\beta(z) = \begin{cases}
-\sqrt{-z} + O(1), & \beta = \infty,\\
\beta + o(1), & \beta\in\R,
\end{cases}
\]
show that $\ti{m}(z)$ is uniquely determined by its zeros and poles via its Hadamard product. To complete the proof of (i) it suffices to note that the Weyl function $m_\beta(z)$ uniquely determines the potential $q$ and the constants $l$ and $\beta$ (see \cite{mar}).

(ii) Similarly, the Herglotz function $m_\beta(z)$ is uniquely determined by its poles $\lam_n^\beta$ and
residues $-\gam_n^\beta$.

(iii) Let $\beta\neq \infty$ and $\phi(\lam_n^\beta,1)$ be given.
Clearly $l$ can be read off by the asymptotics of $\lam_n^\beta$. So let us consider $H$ and $\hat{H}$
associated with $q$ and $\hat{q}$, but with the same $l$. Now, following \cite[Chapter 3]{pt} consider the meromorphic function
\be\label{est_fbeta}
\frac{\big(\phi(z,x) -\hat{\phi}(z,x)\big) \big(\psi^\beta(z,x) -\hat{\psi}^\beta(z,x)\big)}{\phi^\beta(z)}.
\ee
The only poles are at the eigenvalues $\lam_n^\beta$ and the residues are given by
\[
\frac{\big(\phi(\lam_n^\beta,x) -\hat{\phi}(\lam_n^\beta,x)\big)^2}{\phi(\lam_n^\beta,1)\dot{\phi}^\beta(\lam_n^\beta)} \geq 0,
\]
where we have used $\phi(\lam_n^\beta,x) = \phi(\lam_n^\beta,1) \psi^\beta(\lam_n^\beta,x)$ together with our assumption
$\phi(\lam_n^\beta,1) = \hat{\phi}(\lam_n^\beta,1)$. Using our estimates \eqref{estphi}, \eqref{estphi'}, and \eqref{estpsi}, we have
\begin{align*}
\big|\phi(z,x) - \hat{\phi}(z,x)\big| &\leq C \left(\frac{x}{1+ |z|^{1/2} x}\right)^{l+1} \E^{|\im(z^{1/2})| x} \int_0^x \frac{y (\ti{q}(y)+\ti{\hat{q}}(y))}{1 +|z|^{1/2} y} dy,\\
\big|\psi^\beta(z,x) - \hat{\psi}^\beta(z,x)\big| &\leq C \left(\frac{1+ |z|^{1/2} x}{x+ |z|^{1/2}x}\right)^{l} \E^{|\im(z^{1/2})| (1-x)} \int_x^1 \frac{y (\ti{q}(y)+\ti{\hat{q}}(y))}{1 +|z|^{1/2} y} dy.
\end{align*}
Further, together with the estimate from the proof of Theorem~\ref{thmev}, we get
\[
\big|\phi^\beta(z)\big|\geq \big|\phi^\beta_l(z)\big|-\big|\phi^\beta(z)-\phi_l^\beta(z)\big| > \frac{1}{2} \big|\phi_l^\beta(z)\big|
\geq  C_l\frac{\E^{|\im (z)^{1/2}|}}{|z|^{\frac{l}{2}}}
\]
for $|z|= r_n = \big(n+(l+1)/2\big)^2\pi^2$ with $n$ sufficiently large. Observing that for every $x\in(0,1)$
\[
\int_x^1 \frac{y (\ti{q}(y)+\ti{\hat{q}}(y))}{1 +|z|^{1/2} y} dy \leq \frac{1}{1 +|z|^{1/2} x}\int_x^1 y (\ti{q}(y)+\ti{\hat{q}}(y)) dy = O(|z|^{-1/2})
\]
we see that the function \eqref{est_fbeta} is $o(r_n^{-1})$ along these circles and thus Lemma~3.2 from \cite{pt}
shows that the residues are all zero. Hence $\phi(\lam_1^\beta,x) = \hat{\phi}(\lam_1^\beta,x)$ implying $q = \hat{q}$.

If $\beta \neq 0$ and the set $\phi'(\lam_n^\beta,1)$ is given, then the proof remains essentially the same.
The only difference is that we now use $\phi(\lam_n^\beta,x) = -\beta^{-1} \phi'(\lam_n^\beta,1) \psi^\beta(\lam_n^\beta,x)$
if $\beta\ne\infty$. In the case $\beta=\infty$ we need to use $\phi^\beta(z)=\phi'(z,1)$, $\phi(\lam_n^\beta,x) = \phi'(\lam_n^\beta,1) \psi^\beta(\lam_n^\beta,x)$,
and the estimate \eqref{estpsi_B} instead of \eqref{estpsi}.
\end{proof}

\appendix

\section{Some estimates for the spherical Schr\"odinger equation}

In this appendix we want to provide some estimates for the solutions of the spherical Schr\"odinger equation
which are crucial for the main body of our paper. These results are due to Guillot and Ralston \cite{gr}
in the case $l=1$. The analog estimates for arbitrary $l\geq -\frac{1}{2}$ have been stated in \cite{zs}
without proof. However, the estimate (39) from \cite{zs} is clearly wrong in the case $l= -\frac{1}{2}$.
The case of integer $l$ is given in \cite{se}. Since the proof in \cite{se} (and \cite{gr}) uses the explicit
representation of spherical Bessel functions in terms of trigonometric functions we have decided
to provide the details for the general case in this appendix. We assume that $x \in (0, \, 1)$.

\begin{lemma} \label{LaA1}
For $l>-\frac{1}{2}$ the following estimates hold:
\begin{align}\label{estphil}
|\phi_l(z,x)| &\leq C \left(\frac{x}{1+ |z|^{1/2} x}\right)^{l+1} \E^{|\im(z^{1/2})| x},\\ \label{estGl}
|G_l(z,x,y)| &\leq C \left(\frac{x}{1+ |z|^{1/2} x}\right)^{l+1} \left(\frac{1+ |z|^{1/2} y}{y}\right)^l \E^{|\im(z^{1/2})| (x-y)}, \quad y \leq x.
\end{align}
For the case $l=-\frac{1}{2}$  formula \eqref{estphil} remains valid and one has to replace \eqref{estGl} by
\begin{align}\label{a8}
|G_{-1/2}(z,x,y)| \leq &C \left(\frac{xy}{(1+ |z|^{1/2} x)(1+ |z|^{1/2} y)}\right)^{1/2}
\\ \nn & \times
\E^{|\im(z^{1/2})| (x-y)}(1-\log( y)), \quad y \leq x.
\end{align}
\end{lemma}

\begin{proof}
First of all recall \eqref{defphil}, \eqref{defthetal}, and \eqref{a2}, and note that for all $l\ge -\frac{1}{2}$
\begin{align}
\phi_l(z,x) &= \sqrt{\frac{\pi x}{2}} z^{-\frac{2l+1}{4}} J_{l+\frac{1}{2}}(\sqrt{z} x), \qquad \\ \nn
G_l(z,x,y) &= -\frac{\pi}{2} \sqrt{x y} \big(  J_{l+\frac{1}{2}}(\sqrt{z} x) Y_{l+\frac{1}{2}}(\sqrt{z} y) - J_{l+\frac{1}{2}}(\sqrt{z} y) Y_{l+\frac{1}{2}}(\sqrt{z} x)\big)\\
\label{a10}
&= -\frac{\pi \I}{4} \sqrt{x y} \big(  H_{l+\frac{1}{2}}^{(1)}(\sqrt{z} x) H_{l+\frac{1}{2}}^{(2)}(\sqrt{z} y) - H_{l+\frac{1}{2}}^{(1)}(\sqrt{z} y)
H_{l+\frac{1}{2}}^{(2)}(\sqrt{z} x)\big).
\end{align}
Here $H_{l+1/2}^{(1)}(z) = J_{l+1/2}(z) + \I Y_{l+1/2}(z)$ and $H_{l+1/2}^{(2)}(z) = J_{l+1/2}(z) - \I Y_{l+1/2}(z)$
are the Hankel functions of the first and second kind, respectively. Moreover,
recall $Y_{l+1/2}(z) = \sin((l+1/2)\pi)^{-1}(\cos((l+1/2)\pi) J_{l+1/2}(z) - J_{-l-1/2}(z))$ (\cite[(9.1.2)]{as}).

The claimed estimates \eqref{estphil} and \eqref{estGl} can be shown by combining the following asymptotic expansions which yield estimates for
$|z| \ge 1$ and $|z|\leq 1$, respectively. In fact, we have \cite[(9.2.7)-(9.2.10)]{as}
\begin{align}\label{a10'}
H_{l+\frac{1}{2}}^{(1)}(z) &= \sqrt{\frac{2}{\pi z}} \E^{\I(z -(l+1) \pi/2)} \big( 1 + O(1/z) \big), \quad -\pi < \arg(z) < 2\pi,\\
\label{a11}
H_{l+\frac{1}{2}}^{(2)}(z) &= \sqrt{\frac{2}{\pi z}} \E^{-\I(z -(l+1) \pi/2)} \big( 1 + O(1/z) \big),\quad -2\pi < \arg(z) < \pi,
\end{align}
for $z \to \infty$ with the error uniform strictly inside the indicated sectors (\cite[(7.2.1)-(7.2.2)]{wa}). Using the power series for the
Bessel and Neumann functions \cite[(9.1.10)--(9.1.11)]{as} we get
\begin{align}\label{estjl}
J_{l+\frac{1}{2}}(z) &= \frac{1}{\Gamma({l+\frac{3}{2}})} \left(\frac{z}{2}\right)^{l+\frac{1}{2}} (1 + O(z)),\\
Y_{l+\frac{1}{2}}(z) &= \begin{cases}  -\frac{\Gamma({l+\frac{1}{2}})}{\pi} \left(\frac{z}{2}\right)^{-l-\frac{1}{2}} (1 + O(z^{\min(1,2l+1)})), & {l+\frac{1}{2}} >0,\\
\frac{2}{\pi} \log(z) + O(1), & l=-\frac{1}{2}, \end{cases}\label{estyl}
\end{align}
for $z \to 0$.
Thus, \eqref{estphil} immediately follows from \eqref{a10'}, \eqref{a11}, and \eqref{estjl}. The estimate \eqref{estGl} can be also deduced from \eqref{a10'}--\eqref{estyl}.
For instance, using \eqref{estjl}, \eqref{estyl}, and the fact that the function $x\mapsto x(1+x)^{-1}$ is increasing,
we get for $|\eta|\le |\xi| \le 1$
\begin{align*}
& \left| J_{l+\frac{1}{2}}(\xi) Y_{l+\frac{1}{2}}(\eta) - J_{l+\frac{1}{2}}(\eta) Y_{l+\frac{1}{2}}(\xi) \right|\\
&\qquad  \le C \left( \left(\frac{|\xi|}{1+|\xi|}\right)^{l+\frac{1}{2}} \left(\frac{|\eta|}{1+|\eta|}\right)^{-l-\frac{1}{2}}+ \left(\frac{|\eta|}{1+|\eta|}\right)^{l+\frac{1}{2}} \left(\frac{|\xi|}{1+|\xi|}\right)^{-l-\frac{1}{2}} \right)\\
& \qquad\qquad
\le 2 C \left(\frac{|\xi|}{1+|\xi|}\right)^{l+\frac{1}{2}} \left(\frac{|\eta|}{1+|\eta|}\right)^{-l-\frac{1}{2}}.
\end{align*}
Similarly one handles the cases $|\eta|\le 1 \le |\xi|$ and $1 \le |\eta|\le |\xi|$ to obtain the desired result
(for the last one use the expression in terms of Hankel functions).

The estimate \eqref{a8} requires some further considerations, and we
also split this case into three subcases: a)  sufficiently large values of  $|z|^{1/2}y$;
 b) the values of  $|z|^{1/2}x$ are sufficiently large and  the values of  $|z|^{1/2}y$
are bounded;
c) $|z|^{1/2}x$ is bounded.
In all the subcases it is assumed that $y\leq x$.

To prove \eqref{a8} for the subcase of sufficiently large values of $|z|^{1/2}y$
 one can use again the equality \eqref{a10}:
\begin{align}\label{a9}
G_{-1/2}(z,x,y)=-\frac{\pi \I}{4} \sqrt{x y} \big(  H_{0}^{(1)}(\sqrt{z} x) H_{0}^{(2)}(\sqrt{z} y) - H_{0}^{(1)}(\sqrt{z} y)
H_{0}^{(2)}(\sqrt{z} x)\big).
\end{align}
Taking into account \eqref{a10'} and   \eqref{a11} we easily derive \eqref{a8} in the subcase a).

Assume now that the values of $|z|^{1/2}y$ are bounded. Note that  (see \eqref{defthetal}, \eqref{a2})
\[
G_{-1/2}(z,x,y)=\phi_{-\frac{1}{2}}(z,x) \theta_{-\frac{1}{2}}(z,y) - \phi_{-\frac{1}{2}}(z,y) \theta_{-\frac{1}{2}}(z,x),
\]
where
\[
\theta_{-\frac{1}{2}}(z,x)=\sqrt{\frac{\pi x}{2}}\left(Y_0(\sqrt{z} x)-\frac{1}{\pi}\log (z)J_0(\sqrt{z} x)\right).
\]
Using the integral representations \cite[(9.1.18) and (9.1.19)]{as}, we have
\begin{align}\label{a12}&
\ta_{-\frac{1}{2}}(z,x)=\pi^{-3/2}\sqrt{2x}
\int_0^{\pi}\cos\big(\sqrt{z}x\cos \ta\big)\big\{\gam+\log (x)+\log(2\sin^2\ta)\big\}d\ta.
\end{align}
It is immediate from \eqref{a12} that
\begin{align}\label{a13}&
\big|\ta_{-\frac{1}{2}}(z,x)\big| \leq C_1\sqrt{x}(1 - \log (x)) \E^{|\im(z^{1/2})| x} .
\end{align}
Further, since the values of $|z|^{1/2}y$
are bounded,  by \eqref{estphil} and \eqref{a13} we derive
\begin{align}\label{a14}&
\big|\phi_{-\frac{1}{2}}(z,x) \ta_{-\frac{1}{2}}(z,y)\big| \leq C_2 \left(\frac{xy}{1+ |z|^{1/2} x}\right)^{1/2}
(1-\log(y))
\E^{|\im(z^{1/2})| x}
\\ \nn & \leq
C_3 \left(\frac{xy}{(1+ |z|^{1/2} x)(1+ |z|^{1/2} y)}\right)^{1/2}
\E^{|\im(z^{1/2})| (x-y)}(1-\log( y)).
\end{align}

For large values of $|z|^{1/2} x$, Hankel's asymptotic expansions \cite[(9.2.5) and (9.2.6)]{as} yield:
\begin{align}\label{a15}&
\big|\ta_{-\frac{1}{2}}(z,x)\big| \leq C_4 |z|^{-1/4}(1 + \log (|z|)) \E^{|\im(z^{1/2})| x} .
\end{align}
Since $|z|^{1/2}x$ is sufficiently large, the values $|z|$ and $\log(|z|)$ are also large. Moreover, $|z|^{1/2}y$
is bounded and hence $\log(|z|)+2 \log (y)$ is bounded from above.
Thus, \eqref{a15} can be rewritten in the form
\begin{align}\label{a16}&
\big|\ta_{-\frac{1}{2}}(z,x)\big| \leq C_5 \left(\frac{x}{1+ |z|^{1/2} x}\right)^{1/2}
 \E^{|\im(z^{1/2})| x} (1 - \log (y)).
\end{align}
It follows from \eqref{estphil} and \eqref{a16}
that
\begin{align}\label{a17}
\big|\ta_{-\frac{1}{2}}(z,x) \phi_{-\frac{1}{2}}(z,y)\big| \leq &C_6 \left(\frac{xy}{(1+ |z|^{1/2} x)(1+ |z|^{1/2} y)}\right)^{1/2}
\\ &\nn \times
\E^{|\im(z^{1/2})| (x-y)}(1-\log( y)).
\end{align}
Inequalities \eqref{a14} and \eqref{a17} yield \eqref{a8} for the subcase b).

Finally, if  $|z|^{1/2}x$ is bounded, according to  \eqref{estphil}  and \eqref{a13}
we have
\begin{align}\label{a18}&
\big|\ta_{-\frac{1}{2}}(z,x) \phi_{-\frac{1}{2}}(z,y)\big| \leq C_7 \sqrt{xy}(1-\log(x))\leq C_7 \sqrt{xy}(1-\log(y)).
\end{align}
Inequalities \eqref{a14} and \eqref{a18} imply \eqref{a8} for the subcase c).
\end{proof}

To consider  $\frac{\partial}{\partial x} \phi_l(z,x) $ and $\frac{\partial}{\partial x} G_l(z,x) $
one should simply use the asymptotic relations   \cite[(9.2.11)-(9.2.16)]{as} on the derivatives
of the special functions, in addition to
the relations \cite[(9.2.5)-(9.2.10)]{as} which were used in the proof of the previous lemma.
In particular, for sufficiently large values of $|z|$ we have
\begin{align}&\label{a18!}
\left|\frac{\partial}{\partial z}J_{l+\frac{1}{2}}(z) \right|
+\left|\frac{\partial}{\partial z}Y_{l+\frac{1}{2}}(z) \right| \leq C_8\sqrt{\frac{2}{\pi z}}
\E^{|\im(z)|},
\end{align}
and, as $z \to \infty$, we get
\begin{align}
\label{a19}&
\frac{\partial}{\partial z} H_{l+\frac{1}{2}}^{(1)}(z) = i\sqrt{\frac{2}{\pi z}} \E^{\I(z -(l+1) \pi/2)} \big( 1 + O(1/z) \big), \quad -\pi < \arg(z) < 2\pi,\\
\label{a20}&
\frac{\partial}{\partial z} H_{l+\frac{1}{2}}^{(2)}(z) = -i\sqrt{\frac{2}{\pi z}} \E^{-\I(z -(l+1) \pi/2)} \big( 1 + O(1/z) \big),\quad -2\pi < \arg(z) < \pi.
\end{align}
\begin{lemma} \label{LaA2}
For $l>-\frac{1}{2}$ the following estimates hold:
\begin{align}\label{estphilp}
\left|\frac{\partial}{\partial x} \phi_l(z,x) \right| &\leq C \left(\frac{x}{1+ |z|^{1/2} x}\right)^l \E^{|\im(z^{1/2})| x},\\ \label{estGlp}
\left|\frac{\partial}{\partial x} G_l(z,x,y)\right| &\leq C \left(\frac{x}{1+ |z|^{1/2} x}\right)^l \left(\frac{1+ |z|^{1/2} y}{y}\right)^l \E^{|\im(z^{1/2})| (x-y)}, \quad y \leq x.
\end{align}
For  $l=-1/2$ formula \eqref{estphilp} remains true
and one has to replace \eqref{estGlp} by
\begin{align}
\left|\frac{\partial}{\partial x}G_{-1/2}(z,x,y)\right| \leq C \left(\frac{y+ |z|^{1/2} xy)}{x+ |z|^{1/2} xy)}\right)^{1/2}
\E^{|\im(z^{1/2})| (x-y)}(1-\log( y)), \quad y \leq x.\label{a21}
\end{align}
\end{lemma}
\begin{proof} The proof is similar to the proof of Lemma \ref{LaA1} and we shall
prove here only formula \eqref{a21}.  For the subcase of sufficiently large values of $|z|^{1/2}y$
we use  the equality \eqref{a9}. Therefore, in view of \eqref{a10'}, \eqref{a11},
\eqref{a19}, and \eqref{a20} we derive
\begin{align}\label{a22}&
|G_{-1/2}(z,x,y)| \leq C_9
\E^{|\im(z^{1/2})| (x-y)}, \quad y \leq x,
\end{align}
and  \eqref{a21} is immediate.

When the values of $|z|^{1/2}y$ are bounded, formulas \eqref{a13} and \eqref{estphilp}
imply the inequality
\begin{align}\label{a23}&
\left|\frac{\partial}{\partial x} \phi_{-\frac{1}{2}}(z,x) \ta_{-\frac{1}{2}}(z,y)\right| \leq C_{10}
\left(\frac{y+ |z|^{1/2} xy}{x+ |z|^{1/2} xy}\right)^{1/2}
\E^{|\im(z^{1/2})| (x-y)}(1-\log( y)),
\end{align}
which we shall use instead of \eqref{a14}.
For the subcase of sufficiently large values of $|z|^{1/2}x$ and bounded values
of  $|z|^{1/2}y$, it follows from \eqref{defthetal},
\eqref{a16}, and  \eqref{a18!}
that
\begin{align}\label{a24}
\left|\frac{\partial}{\partial x}  \ta_{-\frac{1}{2}}(z,x)\right| \leq C_{11}
\left(\frac{1+ |z|^{1/2} x}{x}\right)^{1/2}\E^{|\im(z^{1/2})| x}(1-\log( y)).
\end{align}
Taking into account \eqref{estphil}, \eqref{a24}, and boundedness of
$|z|^{1/2}y$ we have
\begin{align}\label{a25}
\left|\frac{\partial}{\partial x} \ta_{-\frac{1}{2}}(z,x) \phi_{-\frac{1}{2}}(z,y)\right| \leq  C_{12}
\left(\frac{y+ |z|^{1/2} xy}{x+ |z|^{1/2} xy)}\right)^{1/2}
\E^{|\im(z^{1/2})| (x-y)}(1-\log(y)).
\end{align}
Inequalities \eqref{a23} and \eqref{a25} yield \eqref{a21}
for the subcase of sufficiently large values of $|z|^{1/2}x$ and bounded values
of  $|z|^{1/2}y$.

Finally, if the values of $|z|^{1/2}x$ are bounded, formula \eqref{a12} implies
\begin{align}\label{a26}&
\left|\frac{\partial}{\partial x}  \ta_{-\frac{1}{2}}(z,x)\right| \leq C_{13}x^{-1/2}(1-\log(x)),
\end{align}
and \eqref{a25} follows from \eqref{estphil} and \eqref{a26}.
By \eqref{a23} and \eqref{a25} the inequality \eqref{a21} holds again.
\end{proof}

\noindent
{\bf Acknowledgments.}
We thank Fritz Gesztesy and James Ralston for several helpful discussions. A.K. acknowledges the hospitality and financial
support of the Erwin Schr\"odinger Institute and financial support from the IRCSET PostDoctoral Fellowship Program.

\end{document}